\documentclass[11pt]{article}
\usepackage{amsfonts, amsmath, amssymb, latexsym}

\numberwithin{equation}{section}

\newtheorem{definition}[equation]{Definition}
\newtheorem{theorem}[equation]{Theorem}
\newtheorem{lemma}[equation]{Lemma}

\newtheorem{remark}[equation]{Remark}
\newtheorem{example}[equation]{Example}

\def\C{\mathbb C}
\def\K{\mathbb K} 
\def\Z{\mathbb Z}
\def\K{\mathbb F}
\def\fld{\mathbb K} 

\newcommand{\uqhat}{U_q(\widehat{\mathfrak{sl}}_2)}

\begin{document}

\title{ \bf Irreducible Modules for \\
the Quantum Affine Algebra
 $\uqhat$ and \break  its
Borel Subalgebra $\uqhat^{\geq 0}$} 
\author {Georgia Benkart\footnote{Support from NSF grant \#{}DMS--0245082 is gratefully acknowledged. \hfil \break 
 }  \ and Paul Terwilliger}
\maketitle
\begin{abstract} 
Let
$\uqhat^{\geq 0}$ denote the
Borel subalgebra of the quantum affine algebra
 $\uqhat$.
 We show that the following hold
for any choice of scalars
$\varepsilon_0, \varepsilon_1$  
{f}rom the set $\lbrace 1,-1\rbrace$.
\begin{enumerate}
\item[{\rm (i)}]  Let $V$ be a
finite-dimensional irreducible $\uqhat^{\geq 0}$-module of type
$(\varepsilon_0,\varepsilon_1)$. Then
 the action of 
 $\uqhat^{\geq 0}$ on $V$ extends uniquely
to an action of 
 $\uqhat$ on $V$.
The resulting 
$\uqhat$-module structure on $V$ is irreducible and
of type  $(\varepsilon_0,\varepsilon_1)$.
\item[{\rm (ii)}] Let $V$ be a finite-dimensional irreducible
$\uqhat$-module of type
$(\varepsilon_0,\varepsilon_1)$.  When the
$\uqhat$-action is restricted to 
$\uqhat^{\geq 0}$,
the resulting $\uqhat^{\geq 0}$-module structure on $V$
is irreducible and of type
$(\varepsilon_0,\varepsilon_1)$.
\end{enumerate}
 \end{abstract}

\section{The quantum affine algebra 
$\uqhat$}

\medskip

 The affine Kac-Moody Lie algebra $\widehat {\mathfrak{sl}}_2$ 
  has played an essential  role in diverse areas of mathematics and physics.  
  Elements of  $\widehat {\mathfrak{sl}}_2$  can be represented as  vertex operators,
  which are certain generating functions that appear in the dual resonance models of particle physics  (see \cite{LW} and \cite{FK}).  The algebra $\widehat {\mathfrak{sl}}_2$ also features prominently in the study of Knizhnik-Zamolodchikov equations and conformal field theory (see for example, \cite{cp1} and \cite{EFK}).   Our main object of interest  is a $q$-analogue of $\widehat {\mathfrak{sl}}_2$, the quantum affine algebra $\uqhat$, which also has a representation
 by vertex operators  \cite{FJ} and has many important  connections
with quantum field theory and symmetric functions, in particular with Kostka-Foulkes polynomials
(\cite{NY}, \cite{FKLMM}).   In this paper,  we focus on the 
finite-dimensional irreducible modules of $\uqhat$. 
These modules have been
classified up to isomorphism by V. Chari and A. Pressley \cite{cp}. 
Our aim here is to relate them to 
the finite-dimensional irreducible modules  
of the Borel subalgebra $\uqhat^{\geq 0}$.  
 
 \medskip Throughout the paper   $\K$ will  denote an algebraically closed field. We fix a nonzero scalar
 $q\in \K$ that is not
 a root of unity and adopt  the following notation:

\begin{equation}
\lbrack n \rbrack = {{q^n-q^{-n}}\over {q-q^{-1}}},
\qquad \qquad n=0,1,\ldots
\label{eq:brackndef}
\end{equation}

\begin{definition} 
\label{def:uq}
The quantum affine algebra 
$\uqhat$ is the unital associative $\K$-algebra 
with
generators $e^{\pm}_i$, $K_i^{{\pm}1}$, $i\in \lbrace 0,1\rbrace $
which satisfy  the following relations:
\begin{gather}
K_i K^{-1}_i =
K^{-1}_i K_i =  1,
\label{eq:buq1}
\\
K_0K_1= K_1K_0,
\label{eq:buq2}
\\
K_ie^{\pm}_iK^{-1}_i = q^{{\pm}2}e^{\pm}_i,
\label{eq:buq3}
\\
K_ie^{\pm}_jK^{-1}_i = q^{{\mp}2}e^{\pm}_j, \qquad i\not=j,
\label{eq:buq4}
\\
\lbrack e^+_i, e^-_i\rbrack = {{K_i-K^{-1}_i}\over {q-q^{-1}}},
\label{eq:buq5}
\\
\lbrack e^{\pm}_0, e^{\mp}_1\rbrack = 0,
\label{eq:buq6} \\
(e^{\pm}_i)^3e^{\pm}_j - 
\lbrack 3 \rbrack (e^{\pm}_i)^2e^{\pm}_j e^{\pm}_i 
+\lbrack 3 \rbrack e^{\pm}_ie^{\pm}_j (e^{\pm}_i)^2 - 
e^{\pm}_j (e^{\pm}_i)^3 =0, \qquad i\not=j.
\label{eq:buq7}
\end{gather}
\end{definition}

We call $e^{\pm}_i$, $K_i^{{\pm}1}$, $i\in \lbrace 0,1\rbrace $
the {\it Chevalley generators} for
$\uqhat$ and refer to 
(\ref{eq:buq7})
as the
{\it $q$-Serre relations}.  We denote by
$\uqhat^{\geq 0}$ 
 the subalgebra of $\uqhat$ generated by the elements
 $e_i^+$, $K_i^{\pm 1}$, $i \in \{0,1\}$. We call
$\uqhat^{\geq 0}$ 
 the {\em Borel subalgebra} of $\uqhat$,   
because of its similarity to the universal enveloping algebra 
of the standard  Borel subalgebra of a finite-dimensional simple 
 Lie algebra over the complex numbers. 

\medskip It is apparent from the definitions 
that in  
 $\uqhat$ or $\uqhat^{\geq 0}$, $K_0K_1$ is central
and so by Schur's Lemma must act as a scalar on any finite-dimensional
irreducible module.  Therefore, 
finite-dimensional
irreducible modules of the Borel subalgebra  
are closely related to  the
finite-dimensional
irreducible modules  of the following algebra.

\begin{definition}
\label{def:nn}
\rm
The algebra $U^{\geq 0}$ is the unital
associative $\K$-algebra with generators $R,L, K^{\pm 1}$,
which satisfy the defining relations: 
\begin{eqnarray}
&& \qquad KK^{-1} = 
K^{-1}K =  1,
\label{eq:nnbuq1}
\\
&&\qquad KRK^{-1} = q^{2}R,
\label{eq:nnbuq3}
\\
&&\qquad KLK^{-1} = q^{-2}L,
\label{eq:nnbuq4}
\\
&&R^3L-  
\lbrack 3 \rbrack R^2LR 
+\lbrack 3 \rbrack  RLR^2 - 
LR^3 =0,
\label{eq:nnbuq5}
\\
&&L^3R-  
\lbrack 3 \rbrack  L^2RL 
+\lbrack 3 \rbrack LRL^2 - 
RL^3 =0.
\label{eq:nnbuq7}
\end{eqnarray}
\end{definition}
  
Our first goal is to explain the exact relationship between
the finite-dimensional irreducible
$\uqhat$-modules
and
the finite-dimensional irreducible 
$U^{\geq 0}$-modules.
In order to 
state our results precisely, it is necessary
to make a few comments.

\medskip
 
Let $V$ denote a
finite-dimensional irreducible
$\uqhat$-module.
Then the actions of $K_0$ and
$K_1$ on $V$ are semisimple \cite[Prop. 3.2]{cp}.
Furthermore (also by \cite[Prop. 3.2]{cp}),  there exists an integer
$d \geq 0$ and scalars
$\varepsilon_0,
 \varepsilon_1$ chosen from 
$\lbrace 1,-1\rbrace $
such that
\begin{enumerate}
\item[{\rm (i)}]
the set of distinct eigenvalues of $K_0$  on $V$ is 
$\lbrace \varepsilon_0 q^{2i-d} \mid 0 \leq i \leq d\rbrace $; and 
\item[{\rm (ii)}] 
$K_0K_1-\varepsilon_0 \varepsilon_1 I$ vanishes on $V$.
\end{enumerate}
We call the ordered pair $(\varepsilon_0, \varepsilon_1)$ 
the {\it type} of
 $V$.

\medskip 
Now let  $V$ denote a
finite-dimensional irreducible
$U^{\geq 0}$-module.  As we will see in Section 2, 
the action of $K$ on $V$ is semisimple.  Moreover, 
there exists an integer $d\geq 0$
and a nonzero scalar $\alpha  \in \K$ such that
the set of distinct eigenvalues of $K$ on $V$ is 
$\lbrace \alpha q^{2i-d} \mid 0 \leq i \leq d\rbrace $.
We refer to $\alpha$ as the {\it type} of
 $V$.

\medskip
Our main results concerning
$\uqhat$ and $U^{\geq 0}$ 
are contained in the following two theorems.

\begin{theorem}
\label{thm:nnmain}
Let $V$ denote a finite-dimensional irreducible 
$U^{\geq 0}$-module of type $\alpha$.
Assume $\varepsilon_0, \varepsilon_1$ are scalars in $\lbrace 1,-1\rbrace $.  
Then there exists a unique   
$\uqhat$-module structure on $V$ 
such that  the operators
$e^+_0-R, \ e^+_1-L, \ K^{\pm 1}_0-
\varepsilon_0\alpha^{\mp 1}K^{\pm 1}$, and
 $K^{\pm 1}_1
-\varepsilon_1\alpha^{\pm 1} K^{\mp 1}$ vanish  on $V$.
This $\uqhat$-module structure is irreducible and of
type
$(\varepsilon_0, \varepsilon_1)$.
\end{theorem}

\begin{theorem}
\label{thm:nnmain2}
Let $V$ be a finite-dimensional irreducible
$\uqhat$-module,  
and assume $(\varepsilon_0, \varepsilon_1)$ is its type.  
Let $\alpha $ denote a nonzero scalar in  $\K$.
Then there exists a unique $U^{\geq 0}$-module structure on
$V$ such that the operators 
$e^+_0-R, \ e^+_1-L, \ K^{\pm 1}_0-
\varepsilon_0\alpha^{\mp 1}K^{\pm 1}$, and
$K^{\pm 1}_1
-\varepsilon_1\alpha^{\pm 1} K^{\mp 1}$ vanish on $V$.
This 
$U^{\geq 0}$-module structure is irreducible and of type $\alpha$.
\end{theorem}

\begin{remark}
\rm
Let $\alpha$ be a nonzero scalar in $\K$,
and let $\varepsilon_0, \varepsilon_1$ denote scalars
in $\lbrace 1,-1\rbrace $.
Combining Theorem
\ref{thm:nnmain} and
Theorem
\ref{thm:nnmain2},  we obtain a bijection between
the following two sets: 
\begin{enumerate}
\item[{\rm (i)}]  the isomorphism classes of 
finite-dimensional irreducible $U^{\geq 0}$-modules of type $\alpha$;
\item[{\rm (ii)}]   the isomorphism classes of 
finite-dimensional irreducible
$\uqhat$-modules of type
$(\varepsilon_0, \varepsilon_1)$.
\end{enumerate}
\end{remark}

\medskip
 
\begin{remark}
\rm
As V. Chari and A. Pressley
\cite{cp} have shown,
each finite-dimensional irreducible 
$\uqhat$-module
has a realization as  a tensor product of evaluation modules. 
In our proofs below, we never have occasion to 
invoke this realization.
 In fact, our arguments are quite
elementary
and require only linear algebra.
  It follows from our work and the results of \cite{cp} 
that all of the finite-dimensional irreducible modules for 
$U^{\geq 0}$ and for $\uqhat^{\geq 0}$
can be obtained from tensor products of evaluation modules of $\uqhat$.    
\end{remark}

The plan for the paper is as follows.
In Section 2,  we state some preliminaries concerning
$\uqhat$-modules and
$U^{\geq 0}$-modules.
Sections 3--11 are devoted to
proving 
Theorem \ref{thm:nnmain}.  
In Section 12,  we prove
Theorem \ref{thm:nnmain2}, and in Section 13, we  discuss
irreducible modules for the Borel subalgebra  $\uqhat^{\geq 0}$.

The proof of Theorem \ref{thm:nnmain}  is an adaptation of
a construction
which T. Ito and the second author used to
get $\uqhat$-actions from
a certain type of tridiagonal pair \cite{uq}. Indeed,  the original
motivation for our work came from the study
of tridiagonal pairs 
(\cite{TD00},
\cite{shape})
and the closely related Leonard pairs
(\cite{qSerre},
\cite{LS24},
\cite{lsint},
\cite{TLT:array}, 
\cite{qrac}, 
\cite{aw},
\cite{LS99}).
A Leonard pair is a pair of semisimple 
linear transformations on
a finite-dimensional vector space, each of which acts
tridiagonally on an eigenbasis for the other
\cite[Defn.~1.1]{LS99}.
There is a close connection between Leonard pairs
and 
the orthogonal polynomials that make up the 
terminating branch of the 
Askey scheme
(\cite{KoeSwa},
\cite{TLT:array}, 
\cite[Appendix A]{LS99}).
A tridiagonal pair is a mild generalization of a Leonard pair 
\cite[Defn.~1.1]{TD00}.
See 
\cite{ch},
\cite{ch2}
for related topics.

\section{Preliminaries}
\medskip 

In this section,  we present some background material
on 
irreducible modules for 
$\uqhat$ and 
$U^{\geq 0}$. 
Towards this purpose, we adopt the following conventions.
Assume $V$ is a nonzero finite-dimensional vector space over $\K$.
Let $d$ denote a nonnegative integer.
By a {\it decomposition of $V$ of diameter $d$}, we mean
a sequence $U_0, U_1, \ldots, U_d$ of nonzero
subspaces of $V$ such that
\begin{eqnarray*}
\label{eq:dec}
V = U_0\oplus  U_1 \oplus \cdots \oplus U_d.  
\end{eqnarray*}
Note we do not assume that the spaces  $U_0, U_1, \ldots, U_d$ have  dimension 1.
 For notational convenience we set $U_{-1}:=0$ and $U_{d+1}:=0$.

\medskip 
\begin{lemma}\label{lem:nnweight}
{\rm \cite[Prop.~3.2]{cp}} \ 
Let $V$ denote a finite-dimensional irreducible 
$\uqhat$-module. 
Then there exist
scalars
$\varepsilon_0,
\varepsilon_1$ in 
$\lbrace 1,-1\rbrace $
and   a decomposition
$U_0, \ldots, U_d$ of $V$ 
such that  
\begin{eqnarray}\label{eq:nnkmove}
(K_0-\varepsilon_0q^{2i-d}I)U_i=0  \ \ \hbox{\rm and} \ \ 
 (K_1-\varepsilon_1q^{d-2i}I)U_i=0
\end{eqnarray}
hold for all $i = 0,1,\dots, d$.  
The sequence 
$\varepsilon_0, \varepsilon_1;
U_0, \ldots, U_d
$
is unique. Moreover,  for $0 \leq i \leq d$ we have 
\begin{eqnarray}
&&
e_0^+
U_i \subseteq U_{i+1}, \qquad
e_1^-
U_i \subseteq U_{i+1}, 
\label{eq:altemove1}
\\
&&
e_0^-
U_i \subseteq U_{i-1}, \qquad
e_1^+
U_i \subseteq U_{i-1}.
\label{eq:altemove2}
\end{eqnarray}
\end{lemma}

\begin{remark} \rm  If $\K$ has characteristic 2,  then in Lemma \ref{lem:nnweight}
we view $\lbrace 1,-1\rbrace $ as having a single element.
\end{remark}

\begin{definition}
\label{def:type}
\rm   
 The ordered pair 
 $(\varepsilon_0, \varepsilon_1)$ in Lemma  
\ref{lem:nnweight} is
 the {\it type}
 of $V$ and   $d$ is the {\it diameter} of $V$.
The sequence 
$U_0$,  $\ldots$, $U_d$
is the 
{\it weight space} decomposition
of
 $V$ (relative to $K_0$ and $K_1$). 
\end{definition} 

\begin{lemma}
\label{lem:autuq}
{\rm \cite[Prop.~3.3]{cp}}
For any choice of scalars $\varepsilon_0, \varepsilon_1$ 
{f}rom $\lbrace 1,-1\rbrace $,
there exists an $\K$-algebra
automorphism of 
$\uqhat$ such that
\begin{eqnarray*}
K_i \rightarrow \varepsilon_i K_i,
\qquad \qquad e^+_i\rightarrow e^+_i,
\qquad \qquad e^-_i\rightarrow \varepsilon_i e^-_i
\end{eqnarray*}
for $i \in \lbrace 0,1\rbrace $.
\end{lemma}

\begin{remark} 
\rm
Given a finite-dimensional irreducible
$\uqhat$-module, 
we can alter its type to any other type
by  applying an automorphism as in Lemma  \ref{lem:autuq}.
\end{remark}

\medskip  The next lemma is reminiscent of
Lemma 
\ref{lem:nnweight}.  

\begin{lemma}
\label{lem:s1}
Let $V$ be a  finite-dimensional
irreducible
$U^{\geq 0}$-module.  Then there exist a  nonzero 
$\alpha \in \K$ and  a decomposition
$U_0,  \ldots, U_d$ of $V$ such that
\begin{eqnarray}
\label{eq:s1dec}
(K-\alpha q^{2i-d}I)U_i=0 \qquad (0 \leq i \leq d).
\end{eqnarray}
The sequence $\alpha; U_0,  \ldots, U_d$ is 
unique. Moreover
\begin{eqnarray}
\label{eq:rmove}
RU_i &\subseteq & U_{i+1} \qquad \qquad (0 \leq i \leq d),
\\ 
LU_i &\subseteq & U_{i-1} \qquad \qquad (0 \leq i \leq d).
\label{eq:lmove}
\end{eqnarray}

\end{lemma}
\noindent {\it Proof:}
For $\theta \in \K$, let
$V^{(\theta)} =\{ v \in V \mid Kv=\theta v\}$.  
Using 
(\ref{eq:nnbuq3}), 
(\ref{eq:nnbuq4}) we find that  \begin{eqnarray}
&&
R V^{(\theta)}  \subseteq V^{(q^2 \theta)}, \qquad \qquad
L V^{(\theta)}   \subseteq V^{(q^{-2} \theta)}.
\label{eq:emove}
\end{eqnarray}
Since $\K$ 
is algebraically closed,  $K$ has an eigenvalue in $\K$,
so  $V^{(\theta)} \neq 0$ for
some $\theta \in\K$. Observe 
$\theta\not=0$ since $K$ is invertible on $V$.
The scalars
$\theta, q^{-2}\theta, q^{-4}\theta, \ldots$ are mutually distinct since
$q$ is not a root of unity, and not all of them can be eigenvalues of $K$ on $V$.
Consequently,  there  is
a nonzero $\zeta \in \K$
such that
$V^{(\zeta)} \neq 0$ and 
$V^{(q^{-2}\zeta)} =0$.
There exists an integer $d \geq 0$  such that
$V^{(q^{2i}\zeta)}$ is nonzero
for $0\leq i \leq d$ and zero for $i=d+1$.
Set
$U_i=V^{(q^{2i}\zeta)}$ for $0 \leq i \leq d$ and
define $U_{-1}= 0$, $U_{d+1}=0$.  Note that
\begin{eqnarray}
(K-q^{2i}\zeta I)U_i=0
\qquad \qquad 
(0 \leq i \leq d).
\label{eq:nnkk}
\end{eqnarray}
Line (\ref{eq:s1dec}) is an easy consequence of 
(\ref{eq:nnkk}) by taking
$\alpha :=\zeta q^d$.  Observe that 
$\alpha\not=0$.
Equations (\ref{eq:rmove}) and
(\ref{eq:lmove}) follow from  
(\ref{eq:emove}). 
We claim 
$U_0, \ldots, U_d$ is a decomposition of $V$.
{F}rom the construction,
each of the spaces 
$U_0, \ldots, U_d$ is nonzero.
Lines (\ref{eq:s1dec})--(\ref{eq:lmove})
imply  $\sum_{i=0}^d U_i$ is invariant under  
 $R,L$,$K^{\pm 1}$.   Since $\sum_{i=0}^d U_i$ is nonzero
and $V$ is an irreducible $U^{\geq 0}$-module, it must be that 
 $V=\sum_{i=0}^d U_i$.  
The sum is direct,  since  
$U_0, \ldots, U_d$
are  eigenspaces for $K$ corresponding to distinct eigenvalues. 
Therefore,  $U_0, \ldots, U_d$ is a decomposition of $V$. 
It is 
clear that
the sequence
$\alpha;
U_0,\ldots, U_d
$
is unique.
\hfill $\Box $  

\begin{definition}
\label{def:type2}   \rm In Lemma  
\ref{lem:s1},  $\alpha$  
is said to be  the {\it type} and $d$ the
 {\it diameter} of $V$.
The sequence 
$U_0, \ldots, U_d$ is 
the
{\it weight space} decomposition of
 $V$ (relative to $K$). 
\end{definition} 

\begin{lemma}
\label{lem:autuq2}
For each nonzero $\alpha \in \K$, 
there exists an $\K$-algebra
automorphism of 
$U^{\geq 0}$ such that
\begin{eqnarray*}
K \rightarrow \alpha K,
\qquad \qquad R\rightarrow R,
\qquad \qquad L\rightarrow L. 
\end{eqnarray*}
\end{lemma}
\noindent {\it Proof:}  This is immediate from
Definition
\ref{def:nn}.
\hfill $\Box $  

\begin{remark}
\rm   
Given a 
 finite-dimensional irreducible
$U^{\geq 0}$-module,
we can change its 
 type to any other type
by  applying an automorphism from 
Lemma
\ref{lem:autuq2}.
\end{remark}

\section{An outline of the proof for Theorem 
\ref{thm:nnmain}}

Our proof of 
Theorem \ref{thm:nnmain} will consume  most of the paper 
from Section 4 to Section 11.
Here we sketch an overview of the argument.

\medskip
 Let $V$ denote a finite-dimensional irreducible
$U^{\geq 0}$-module of type $\alpha$.
For any choice of  $\varepsilon_0,
 \varepsilon_1$  {f}rom 
$\lbrace 1,-1\rbrace $, we begin the construction of 
the $\uqhat$-action 
on $V$ by requiring that the operators 
$e^+_0-R$, $e^+_1-L$, $K^{\pm 1}_0-\varepsilon_0
\alpha^{\mp 1}K^{\pm 1}$,
 $K^{\pm 1}_1
-\varepsilon_1\alpha^{\pm 1} K^{\mp 1}$ vanish  on $V$.
This gives the actions of the elements 
$e^+_0, e^+_1, K^{\pm 1}_0,
 K^{\pm 1}_1$ on $V$.
We define the actions of $e^-_0, e^-_1$ on
$V$ as follows. \ \ 
First we prove  that 
$K+R$ and $K^{-1}+L$ act semisimply on $V$.
Then we show that the set of distinct eigenvalues of
$K+R$ (resp. $K^{-1}+L$) on $V$
is $\lbrace \alpha q^{2i-d}\mid 0 \leq i \leq d\rbrace$
(resp. $\lbrace \alpha^{-1} q^{d-2i}\mid 0 \leq i \leq d\rbrace$), where
$d$ is the diameter of $V$.
For $0 \leq i \leq d$,  we let
$V_i$ (resp. $V^*_i$)
denote the eigenspace of $K+R$  (resp.  of $K^{-1}+L$)
on $V$ associated with the eigenvalue
$\alpha q^{2i-d}$  (resp.
$\alpha^{-1} q^{d-2i}$).
Then $V_0, \ldots, V_d$ \ (resp. $V^*_0, \ldots, V^*_d$) is
a 
decomposition of $V$.   To motivate what comes next, we mention
that  the weight space decomposition $U_0, \ldots, U_d$ of $V$
satisfies
\begin{eqnarray*}
U_i = (V^*_0+\cdots + V^*_i) \cap (V_i + \cdots + V_{d})
\qquad \qquad (0 \leq i \leq d).
\end{eqnarray*}
For $0 \leq i \leq d$,  we define
\begin{eqnarray*}
W_i &=& (V^*_0+\cdots + V^*_i) \cap (V_0 + \cdots + V_{d-i}),
\\
W^*_i &=& (V^*_{d-i}+\cdots + V^*_d) \cap (V_i + \cdots + V_{d}).
\end{eqnarray*}
We argue  that both  $W_0, \ldots, W_d$ and 
$W^*_0, \ldots, W^*_d$ are  decompositions of $V$.
Therefore, there exist linear transformations 
$B:V\rightarrow V$ (resp.  
$B^*:V\rightarrow V$)  
such that
for $0 \leq i \leq d$, $W_i $ (resp. $W^*_i$) is an eigenspace
for $B$ (resp. $B^*$) with eigenvalue $q^{2i-d}$ (resp. $q^{d-2i}$).
We let  $e^-_0$ (resp. $e^-_1$)
act on $V$ as
$\alpha^{-1}I-K^{-1}B$
times
$\varepsilon_0 
q^{-1}(q-q^{-1})^{-2}$
(resp.  
$\alpha I- KB^*$ times
$\varepsilon_1
q^{-1}(q-q^{-1})^{-2}$).
We display some relations that are
satisfied by $B,B^*$,  and
the generators of $U^{\geq 0}$.
Using these relations,  we argue that the above actions
of
$e^{\pm}_0$,
$e^{\pm}_1$,
$K_0^{\pm 1}$,
$K_1^{\pm 1}$
satisfy the defining relations for
$\uqhat$.
In this way,  we obtain the required action of
$\uqhat$ on $V$.

\section{The elements $A$ and $A^*$}

\medskip
\noindent 
As we proceed with our investigation of
$U^{\geq 0}$, we find it convenient  to work with the elements
$K+R$ and $K^{-1}+L$  instead  $R $ and $L$. 
Hence we are led to the following definition. 

\begin{definition}
\label{def:aas}
\rm
Let $A$ and $A^*$ denote the following elements of
$U^{\geq 0}$:
\begin{eqnarray}
A=K+R, \qquad \qquad A^*=K^{-1}+L.
\end{eqnarray}
We observe $A, A^*, K^{\pm 1}$ form a generating set
for $U^{\geq 0}$.
\end{definition}

\begin{lemma}
The following relations hold in $U^{\geq 0}$:
\begin{eqnarray}
\frac{qK^{-1}A-q^{-1}AK^{-1}}{q-q^{-1}}&=& 1,
\label{eq:nnkeq1}
\\
\frac{qKA^*-q^{-1}A^*K}{q-q^{-1}}&=& 1.
\label{eq:nnkeq4}
\end{eqnarray}
\end{lemma}
\noindent {\it Proof:}
To verify 
(\ref{eq:nnkeq1}),
substitute $A=K+R$ and simplify the result using
(\ref{eq:nnbuq3}). Relation 
 (\ref{eq:nnkeq4}) can be verified similarly using (\ref{eq:nnbuq4}).
\hfill $\Box $  

\begin{remark}\rm Each of the equations \eqref{eq:nnkeq1} and \eqref{eq:nnkeq4}
is essentially an instance of the $q$-Weyl relation, which is
the defining relation of the $q$-Weyl algebra (see \cite{G},  for example).  
A presentation of
$\uqhat$ which involves the $q$-Weyl relations and the $q$-Serre relations
can be found in \cite{uq}.
\end{remark}
\begin{lemma}
\label{lem:aasqs} The elements $A$, $A^*$ in Definition \ref{def:aas}  satisfy these relations:
\begin{eqnarray}
A^3A^*-\lbrack 3 \rbrack A^2A^*A+\lbrack 3 \rbrack AA^*A^2-A^*A^3&=&0,
\label{eq:nnaserre}
\\
A^{*3}A-\lbrack 3 \rbrack A^{*2}AA^*+\lbrack 3 \rbrack A^*AA^{*2}-AA^{*3}&=&0.
\label{eq:nnasserre}
\end{eqnarray}
\end{lemma}
\noindent {\it Proof:}
To verify  
\eqref{eq:nnaserre},  substitute 
 $A=K+R$ and
 $A^*=K^{-1}+L$,  and simplify the result using
(\ref{eq:nnbuq3})--(\ref{eq:nnbuq5}).  
Line \eqref{eq:nnasserre} can be checked in the same way.

\hfill $\Box $ 

\begin{lemma}
\label{lem:amove}
Let $V$ be a finite-dimensional irreducible 
$U^{\geq 0}$-module with type 
$\alpha $  and weight space decomposition $U_0, \ldots, U_d$.
Then for $0 \leq i \leq d$ we have 
\begin{eqnarray}
\label{eq:aup}
(A-\alpha q^{2i-d}I)U_i &\subseteq & U_{i+1},
\\
\label{eq:asdown}
(A^*-\alpha^{-1} q^{d-2i}I)U_i &\subseteq & U_{i-1}.
\end{eqnarray}
\end{lemma}
\noindent {\it Proof:}  Relation 
(\ref{eq:aup}) follows directly  from
(\ref{eq:s1dec}), (\ref{eq:rmove}),  and the fact that
$A=K+R$.   Line
(\ref{eq:asdown}) is a consequence of
(\ref{eq:s1dec}), (\ref{eq:lmove}),  and
$A^*=K^{-1}+L$.
\hfill $\Box $ 

\begin{lemma}
\label{lem:ass} \
Let $V$ denote a finite-dimensional irreducible 
$U^{\geq 0}$-module of type $\alpha$ and diameter $d$.
Then the elements $A$ and  $A^*$ in Definition \ref{def:aas} act semisimply on $V$.
The set of distinct eigenvalues for $A$ on $V$ is
$\lbrace \alpha q^{2i-d} \mid 0 \leq i \leq d\rbrace $.
 The set of distinct eigenvalues for $A^*$ on $V$ is
$\lbrace \alpha^{-1} q^{d-2i} \mid 0 \leq i \leq d\rbrace $.
\end{lemma}
\noindent {\it Proof:} 
It is apparent from (\ref{eq:aup}) that the product  $\prod_{i=0}^d (A-\alpha q^{2i-d}I)$
vanishes on $V$. Since
the 
scalars $\alpha q^{2i-d}$ $( 0 \leq i \leq d)$
are mutually distinct,  we find $A$ is semisimple on $V$.
It is clear from 
(\ref{eq:aup})  that the complete set of
distinct eigenvalues of $A$ on $V$ is
$\lbrace \alpha q^{2i-d} \mid 0 \leq i \leq d\rbrace $.
Our assertions concerning $A^*$ follow similarly.
\hfill $\Box $  \\

\section{The eigenspace decompositions for $A$ and $A^*$}

\begin{definition}
\label{def:vivis}
\  \rm
Assume $V$ is  a finite-dimensional irreducible 
$U^{\geq 0}$-module.
Referring to 
Lemma
\ref{lem:ass}, 
we let  $V_i$ (resp. $V^*_i$)
denote the eigenspace of $A$ \ (resp. $A^*$)
on $V$ corresponding to the eigenvalue
$\alpha q^{2i-d}$ \ (resp.
$\alpha^{-1} q^{d-2i}$)  for $0 \leq i \leq d$.
Then 
$V_0, \ldots, V_d$ and $V^*_0, \ldots, V^*_d$ are each
decompositions  of
  $V$.
\end{definition}

\begin{lemma}
\label{lem:uvsv}
Let $V$ denote a  finite-dimensional
irreducible 
$U^{\geq 0}$-module with weight space decomposition 
$U_0, \ldots, U_d$.
Suppose $V_0, \ldots, V_d$
and $V^*_0, \ldots, V^*_d$ are the decompositions in Definition
\ref{def:vivis}. Then for
$0 \leq i \leq d$ we have 
\begin{enumerate}
 \item[{\rm (i)}]  $U_i+  \cdots + U_d = V_i + \cdots + V_d$,
 \item[{\rm (ii)}]  $U_0+  \cdots + U_i = V^*_0 + \cdots + V^*_i$,
 \item[{\rm (iii)}] 
$U_i= (V^*_0 + \cdots + V^*_i)\cap (V_i + \cdots + V_d)$. 
\end{enumerate}
\end{lemma}
\noindent {\it Proof:} (i)
Set $X_i = \sum_{j=i}^d U_j$ and
$X'_i = \sum_{j=i}^dV_j$.
We show $X_i = X'_i$.
Assume $\alpha $ is  the type of $V$,  and
let $T_i : =\prod_{j=i}^d (A-\alpha q^{2j-d}I)$.
Then 
$X'_i = \lbrace v \in V \mid  T_iv=0\rbrace$,
and 
$T_i X_i =0$ by
(\ref{eq:aup})
so $X_i \subseteq X'_i$.
Now let  $S_i :=\prod_{j=0}^{i-1} (A-\alpha q^{2j-d}I)$.
Observe that  
$S_i V=X'_i$,
and 
$S_i V \subseteq X_i$ 
by 
(\ref{eq:aup})
so $X'_i \subseteq X_i$.  We conclude that $X_i = X'_i$ holds as required.
 
 (ii)
Mimic  the proof of (i).
 
 (iii) Combine parts (i) and (ii) above.
\hfill $\Box $

\begin{lemma}
\label{lem:kvi}
Let $V$ be a finite-dimensional
irreducible
$U^{\geq 0}$-module of type $\alpha$ and diameter $d$. Assume  $V_0, \ldots, V_d$
and $V^*_0, \ldots, V^*_d$  are the decompositions in 
Definition
\ref{def:vivis}.
Then for $0 \leq i \leq d$ we have 
\begin{enumerate}
 \item[{\rm (i)}]  $(K^{-1}-\alpha^{-1} q^{d-2i}I)V_i \subseteq V_{i+1}$,
 \item[{\rm (ii)}]  $(K-\alpha q^{2i-d}I)V_i \subseteq 
V_{i+1} + \cdots + V_d$,
 \item[{\rm (iii)}] $(K-\alpha q^{2i-d}I)V^*_i \subseteq V^*_{i-1}$,
 \item[{\rm (iv)}]  $(K^{-1}-\alpha^{-1} q^{d-2i}I)V^*_i \subseteq
V^*_0 + \cdots + V^*_{i-1}$.
\end{enumerate}
\end{lemma}
\noindent {\it Proof:} 
(i)
Recall that  
for $0 \leq i \leq d$, $V_i$ is an eigenspace for $A$ with corresponding 
eigenvalue $\alpha q^{2i-d}$.
Therefore it suffices to show that 
\begin{eqnarray}
(A-\alpha q^{2i+2-d}I)(K^{-1}-\alpha^{-1}q^{d-2i}I)
\label{eq:nnwantk}
\end{eqnarray}
vanishes on $V_i$ for $0 \leq i \leq d$. Since
 $A-\alpha q^{2i-d}I$ vanishes on 
$V_i$,  so does 
\begin{eqnarray}
(K^{-1}-\alpha^{-1}q^{d-2i-2}I)(A-\alpha q^{2i-d}I). 
\label{eq:nnbak}
\end{eqnarray}  Using 
(\ref{eq:nnkeq1}) we see that
\begin{eqnarray}
qK^{-1}A-q^{-1}AK^{-1}-(q-q^{-1})I 
\label{eq:nnqabk}
\end{eqnarray}
is 0 on $V_i$.
Subtracting
(\ref{eq:nnbak})
from
$q^{-1}$ times 
(\ref{eq:nnqabk}) 
we find that 
(\ref{eq:nnwantk}) vanishes on $V_i$.
Hence
$(K^{-1}-\alpha^{-1}q^{d-2i}I)V_i \subseteq V_{i+1}$ for  $0 \leq i \leq d$.

  Part (ii) follows from (i) above, while 
  (iii) can be obtained by an argument similar  to the proof of (i).  
  Finally, (iv) is a consequence of (iii).    \hfill $\Box $

\begin{lemma}
\label{lem:avi}
Let $V$ be a  finite-dimensional irreducible
$U^{\geq 0}$-module. Let the 
decompositions $V_0, \ldots, V_d$
and $V^*_0, \ldots, V^*_d$ be as in 
 Definition
\ref{def:vivis}. 
Then  for $0 \leq i \leq d$ we have
\begin{enumerate}
 \item[{\rm (i)}]  $A^*V_i \subseteq V_{i-1} + V_i + V_{i+1}$,
 \item[{\rm (ii)}]  $AV^*_i \subseteq V^*_{i-1} + V^*_i + V^*_{i+1}$.
\end{enumerate}
\end{lemma}

\noindent {\it Proof:} 
(i)
Let $\alpha$ denote the type of $V$. 
Recall that for $0 \leq i \leq d$, $V_i$ is
the eigenspace for $A$ corresponding to the eigenvalue
$\alpha q^{2i-d}$.
Therefore it suffices to show
\begin{eqnarray*}
\label{eq:aaaneed}
\left (A-\alpha q^{2i-2-d}I\right)
\left (A-\alpha q^{2i-d}I\right)
\left (A-\alpha q^{2i+2-d}I\right)
A^*V_i=0
\end{eqnarray*}
for $0 \leq i \leq d$.  
For $v \in V_i$  we have
\begin{eqnarray*}
0 &=&
\left (A^3A^*-\lbrack 3\rbrack A^2A^*A+\lbrack 3\rbrack AA^*A^2-A^*A^3\right)v
\qquad \qquad \quad 
(\mbox{by}\; (\ref{eq:nnaserre}))
\\
&=&
\left(A^3A^*-\lbrack 3\rbrack A^2A^*\alpha q^{2i-d}+
\lbrack 3\rbrack AA^*\alpha^2q^{4i-2d}-A^*\alpha^3q^{6i-3d}\right)v
\\
&=& 
\left (A-\alpha q^{2i-2-d}I\right)
\left (A-\alpha q^{2i-d}I\right)
\left(A-\alpha q^{2i+2-d}I\right)A^*v,
\end{eqnarray*}
which gives the desired result. 

 (ii)  This can be argued analogously.  \hfill $\Box $ \\

\section{Yet two more decompositions}

\begin{definition}
\label{def:wdec}
\rm 
Let $V$ denote a  finite-dimensional irreducible $U^{\geq 0}$-module.
 Assume 
  $V_0, \ldots, V_d$
and $V^*_0, \ldots, V^*_d$ are the decompositions from
Definition \ref{def:vivis}.
For $0 \leq i \leq d$,   we define
\begin{eqnarray}
\label{eq:widef}
W_i &=& (V^*_0+\cdots + V^*_i) \cap (V_0 + \cdots + V_{d-i}),
\\
W^*_i &=& (V^*_{d-i}+\cdots + V^*_d) \cap (V_i + \cdots + V_{d}).
\end{eqnarray}
\end{definition}

\medskip 
 
Our next goal is to show that the spaces 
$W_0, \ldots, W_d$
and $W^*_0, \ldots, W^*_d$ in Definition 
\ref{def:wdec} afford
decompositions of $V$.
Towards this purpose, the following definition will be useful.

\begin{definition}
\label{def:nnVijS99}
\rm
Let $V$ be a finite-dimensional irreducible
$U^{\geq 0}$-module. Assume 
$V_0, \ldots, V_d$ and 
$V^*_0, \ldots, V^*_d$
are the decompositions 
in
Definition \ref{def:vivis}.
Set
\begin{equation}
W(i,j) =  \Biggl(\sum_{h=0}^i V^*_h\Biggr)\bigcap
 \Biggl(\sum_{k=0}^j V_k\Biggr)
\label{eq:nndefofvijS99}
\end{equation}
for all integers $i,j$. We interpret the left sum in
(\ref{eq:nndefofvijS99}) to be 0 (resp. $V$) if $i<0$  (resp. $i>d $).
Likewise the right sum in
(\ref{eq:nndefofvijS99}) is interpreted to be 0 (resp. $V$) if $j<0$  (resp. $j>d$).
\end{definition}

\begin{example}
\label{lem:thevijbasicfactsS99}
With reference to Definition \ref{def:nnVijS99},
the following hold.
\begin{enumerate}
\item[{\rm (1)}] 
$W(i,d) = 
 V^*_0+V^*_1+\cdots+V^*_i \qquad (0 \leq i \leq d)$.
\item[{\rm (2)}] 
$W(d,j) = 
 V_0+V_1+\cdots+V_{j} \qquad \  (0 \leq j \leq d)$.
\end{enumerate}
\end{example}
\noindent {\it Proof:}
To obtain (1),   set $j=d$ in 
(\ref{eq:nndefofvijS99}) and recall that $\sum_{k=0}^d V_k =V$; 
for (2), use $i = d$ and  $\sum_{h=0}^d V_h^* =V$.
\hfill $\Box $

\begin{lemma}
\label{lem:howaactsonvijS99}
Let $V$ denote a finite-dimensional irreducible
$U^{\geq 0}$-module with type $\alpha $ and diameter
$d$.
Then for the spaces $W(i,j)$ in  
Definition
\ref{def:nnVijS99},  we have  \begin{enumerate}
 \item[{\rm (i)}] 
$(A-\alpha q^{2j-d}I)W(i,j) \subseteq W(i+1,j-1)$,
\item[{\rm (ii)}] 
$(A^*-\alpha^{-1}q^{d-2i}I)W(i,j) \subseteq W(i-1,j+1)$,
\item[{\rm (iii)}] 
$(K^{-1}-\alpha^{-1} q^{d-2i}I)W(i,j) \subseteq W(i-1,j+1)$,
\item[{\rm (iv)}] 
$(K-\alpha q^{2i-d}I)W(i,j) \subseteq \sum_{h=1}^i W(i-h,j+h)$,
\end{enumerate}
for  $0 \leq i,j \leq d$.  
\end{lemma}
\noindent {\it Proof:}
(i)  Using 
Lemma
\ref{lem:avi} (ii), 
we find that  
\begin{equation}
\left (A-\alpha q^{2j-d}I\right)\sum_{h=0}^i V^*_h \subseteq 
\sum_{h=0}^{i+1} V^*_h.
\label{eq:nnaminusthetajAS99}
\end{equation} 
Because each $V_k$ is an eigenspace
for $A$ with eigenvalue $\alpha q^{2k-d}$, we have 
\begin{equation}
\left (A-\alpha q^{2j-d}I\right)\sum_{k=0}^j V_k = 
\sum_{k=0}^{j-1} V_k.
\label{eq:nnaminusthetajBS99}
\end{equation} 
Evaluating $(A-\alpha q^{2j-d}I)W(i,j)$ using 
(\ref{eq:nndefofvijS99})--(\ref{eq:nnaminusthetajBS99}), 
 we see that  it is contained in $W(i+1,j-1)$.
 
(ii)  This part can be demonstrated  using  the  relations  
\begin{eqnarray}
\left(A^*-\alpha^{-1}q^{d-2i}I\right)\sum_{h=0}^i V^*_h &=& 
\sum_{h=0}^{i-1} V^*_h, 
\label{eq:nnaminusthetajAS99ag}  \\
\left (A^*-\alpha^{-1}q^{d-2i}I\right)\sum_{k=0}^j V_k &\subseteq& 
\sum_{k=0}^{j+1} V_k 
\label{eq:nnaminusthetajBS99ag}
\end{eqnarray}
in conjunction with (\ref{eq:nndefofvijS99}).
 
(iii) {F}rom
Lemma 
\ref{lem:kvi} (iv) and Lemma \ref{lem:kvi} (i), 
we find that  
\begin{eqnarray}
\left (K^{-1}-\alpha^{-1} q^{d-2i}I\right)\sum_{h=0}^i V^*_h &\subseteq& 
\sum_{h=0}^{i-1} V^*_h,
\label{eq:nnkaminusthetajAS99} \\
\left (K^{-1}-\alpha^{-1} q^{d-2i}I\right)\sum_{k=0}^j V_k &\subseteq&  
\sum_{k=0}^{j+1} V_k.
\label{eq:nnkaminusthetajBS99}
\end{eqnarray}
Evaluating 
 $(K^{-1}-\alpha^{-1} q^{d-2i}I)W(i,j)$
using (\ref{eq:nndefofvijS99}),
(\ref{eq:nnkaminusthetajAS99}), and (\ref{eq:nnkaminusthetajBS99}),  
we see that it is contained in 
$W(i-1,j+1)$.
 
(iv) This assertion follows from (iii).
\hfill $\Box $

\begin{lemma}
\label{lem:nnmanyvijzeroS99}
For the spaces $W(i,j)$ in 
Definition \ref{def:nnVijS99}, 
\begin{equation}
W(i,j)= 0 \quad \hbox{if}\quad i+j<d,\qquad  \qquad (0 \leq i,j\leq d).
\label{eq:nnmanyvijzeroS99}
\end{equation}
\end{lemma}
\noindent {\it Proof:} Lemma 
\ref{lem:howaactsonvijS99} implies that for $0 \leq r < d$ 
 the sum 
\begin{equation}
W:=  W(0,r)+W(1,r-1)+\cdots +W(r,0)
\label{eq:nnmanyvijzeroAS99}
\end{equation}
is invariant under $A, A^*$,  and $K^{\pm 1}$.
Since $A,A^*, K^{\pm 1}$
is a generating set for $U^{\geq 0}$,  we find that 
$W$ is a 
$U^{\geq 0}$-submodule of $V$.
Because $V$ is an irreducible $U^{\geq 0}$-module,  we have
$W=0$ or $W=V$.     By 
(\ref{eq:nndefofvijS99}), 
 each term in 
(\ref{eq:nnmanyvijzeroAS99}) is contained in 
\begin{equation}
V_0+V_1+\cdots +V_r,
\label{eq:nnmanyvijzeroBS99}
\end{equation}
so $W$ is contained in 
(\ref{eq:nnmanyvijzeroBS99}). The containment of $W$ in 
 $V$ is proper 
since 
$r<d$.
Thus   
$W=0$ 
 and 
(\ref{eq:nnmanyvijzeroS99}) follows.
\hfill $\Box $  

\begin{lemma}
\label{lem:wone}
Let $V$ denote a finite-dimensional irreducible 
$U^{\geq 0}$-module. Then 
the sequence 
 $W_0,  \ldots, W_d$ from  Definition
\ref{def:wdec}
 is a decomposition of $V$.
\end{lemma}
\noindent {\it Proof:}
First we argue that $\sum_{i=0}^d W_i=V$.
Comparing
(\ref{eq:widef}) and
(\ref{eq:nndefofvijS99}),  
we find that $W_i=W(i,d-i)$ for $0 \leq i \leq d$; 
by this and
 Lemma
\ref{lem:howaactsonvijS99},  we have that   
$\sum_{i=0}^d W_i$ is invariant under
$A, A^*, K^{\pm 1}$.
Because  $A,A^*, K^{\pm 1}$
generate $U^{\geq 0}$, 
 $\sum_{i=0}^d W_i$  must be $0$
or $V$.   
Observe that 
$\sum_{i=0}^d W_i$
 contains 
 $W_0 =V^*_0 \neq 0$, 
 so
$\sum_{i=0}^d W_i =V$.
To  show that  the sum 
$\sum_{i=0}^d W_i$ 
is direct, we prove that 
\begin{eqnarray*} 
 (W_{0}+W_{1}+\cdots +W_{i-1})\cap W_{i} =0
\end{eqnarray*}
for $1 \leq i \leq d$.
 From the construction,  
\begin{eqnarray*}
W_{j} \subseteq V^*_0+V^*_1+\cdots +V^*_{i-1}
\end{eqnarray*}
for $0 \leq j \leq i-1$,  and
\begin{eqnarray*}
W_i\subseteq V_0+V_1+\cdots +V_{d-i}.
\end{eqnarray*}
Therefore
\begin{eqnarray*} 
&&(W_{0}+W_{1}+\cdots +W_{i-1})\cap W_{i}
\\
&& \qquad \qquad  \qquad \subseteq
(V^*_0+V^*_1+\cdots +V^*_{i-1})\cap
(V_0+V_{1}+\cdots +V_{d-i}) \qquad \qquad 
\\
&& \qquad \qquad  \qquad = W(i-1,d-i)
\\
&& \qquad \qquad \qquad  = 0,
\end{eqnarray*}
in view of Lemma 
\ref{lem:nnmanyvijzeroS99}.
We have now shown that 
the sum 
$\sum_{i=0}^d W_i$ is direct. Next we 
argue that  $W_i \neq 0$ for $0 \leq i \leq d$.
We have 
$W_0 = V_0^* \neq 0$ and 
$W_d=V_0 \neq 0$.   Suppose there exists an integer $i$ $(1 \leq i \leq d-1)$
such that $W_i=0$.
By Lemma
\ref{lem:howaactsonvijS99},  the sum
$\sum_{h=0}^{i-1}W_i$ is a $U^{\geq 0}$-submodule of $V$
since it is invariant under each of $A,A^*,K^{\pm 1}$.  Therefore 
 $\sum_{h=0}^{i-1}W_i=0$  or
 $\sum_{h=0}^{i-1}W_i=V$.
But 
 $\sum_{h=0}^{i-1}W_i\neq 0$ since $i\geq 1$ and $W_0\neq 0$; 
and 
 $\sum_{h=0}^{i-1}W_i\neq V$ since $i-1<d$ and $W_d\neq 0$.
We have reached a contradiction. Consequently,  $W_i\not=0$ for $0 \leq i \leq d$.
Thus,  the sequence $W_0, \ldots, W_d$ is a decomposition
of $V$.
\hfill $\Box $  

\begin{lemma}
\label{lem:wsone}
Let $V$ denote a finite-dimensional irreducible 
$U^{\geq 0}$-module. Then
the sequence 
 $W^*_0,  \ldots, W^*_d$ from Definition
\ref{def:wdec}
 is a decomposition of $V$.
\end{lemma}
\noindent {\it Proof:}  Imitate  the proof of 
Lemma
\ref{lem:wone}.
\hfill $\Box $ \\

We record a few helpful facts for later use.

\begin{lemma}
\label{lem:aw}
Assume  $V$ is  a finite-dimensional
irreducible $U^{\geq 0}$-module of type $\alpha $ and diameter $d$.
Then for the decomposition $W_0, \ldots, W_d$ in 
Definition
\ref{def:wdec},  the following hold for  $0 \leq i\leq d$.
\begin{enumerate}
\item[{\rm (i)}] 
$(A-\alpha q^{d-2i}I)W_i \subseteq W_{i+1}$.
\item[{\rm (ii)}] 
$(A^*-\alpha^{-1} q^{d-2i}I)W_i \subseteq W_{i-1}$.
\item[{\rm (iii)}] 
$(K^{-1}-\alpha^{-1} q^{d-2i}I)W_i \subseteq W_{i-1}$.
\item[{\rm (iv)}] 
$(K-\alpha q^{2i-d}I)W_i \subseteq W_0 + \cdots + W_{i-1}$. 
\end{enumerate}
\end{lemma}
\noindent {\it Proof:} 
Set $j=d-i$ and $W(i,d-i)=W_i$ in Lemma
\ref{lem:howaactsonvijS99}.
\hfill $\Box $

\begin{lemma}
\label{lem:aws}
Let $V$ denote a finite-dimensional irreducible
$U^{\geq 0}$-module of type $\alpha $ and diameter $d$.
Then for the decomposition
$W^*_0, \ldots, W^*_d$ in 
Definition
\ref{def:wdec},  the following  hold for $0 \leq i \leq d$.
\begin{enumerate}
\item[{\rm (i)}] 
$(A-\alpha q^{2i-d}I)W^*_i \subseteq W^*_{i+1}$.
\item[{\rm (ii)}] 
$(A^*-\alpha^{-1} q^{2i-d}I)W^*_i \subseteq W^*_{i-1}$.
\item[{\rm (iii)}]  
$(K-\alpha q^{2i-d}I)W^*_i \subseteq W^*_{i+1}$.
\item[{\rm (iv)}]
$(K^{-1}-\alpha^{-1} q^{d-2i}I)W^*_i \subseteq W^*_{i+1} + \cdots + W^*_d$.
\end{enumerate}
\end{lemma} 

\begin{lemma}
\label{lem:wsum}
Assume $V$ is a  finite-dimensional irreducible
$U^{\geq 0}$-module.
Let the
decompositions
$V_0, \ldots, V_d$ and 
$V^*_0, \ldots, V^*_d$
be as  in Definition
\ref{def:vivis}. Let
the decompositions $W_0, \ldots, W_d$ 
and $W^*_0, \ldots, W^*_d$ be as  in
Definition
\ref{def:wdec}.
Then the following 
hold for   $0 \leq i \leq d$.
\begin{enumerate}
\item[{\rm (i)}]  $W_i + \cdots +W_d = V_0 + \cdots + V_{d-i}$.
\item[{\rm (ii)}]  $W_0 + \cdots +W_i = V^*_0 + \cdots + V^*_i$.
\item[{\rm (iii)}]  $W^*_i + \cdots +W^*_d = V_i + \cdots + V_d$.
\item[{\rm (iv)}]  $W^*_0 + \cdots + W^*_i = V^*_{d-i} + \cdots + V^*_d$.
\end{enumerate}
\end{lemma}
\noindent {\it Proof:}
(i) Set  $X_i := \sum_{j=i}^{d} W_j$ and
$X'_i = \sum_{j=0}^{d-i}V_j$.
We show $X_i=X'_i$.
Let $\alpha$ denote the type of $V$ and set
$T_i: =\prod_{j=0}^{d-i} (A-\alpha q^{2j-d}I)$.
Then $X'_i = \lbrace v \in V \mid T_iv=0\rbrace$, and 
$T_i X_i =0$ by
Lemma
\ref{lem:aw} (i), 
so   $X_i \subseteq X'_i$.
Now let 
$S_i: =\prod_{j=d-i+1}^{d} (A-\alpha q^{2j-d}I)$.
Observe that
$S_i V=X'_i$,
and 
$S_i V \subseteq X_i$ 
by 
Lemma 
\ref{lem:aw} (i), 
so $X'_i \subseteq X_i$. We conclude  $X_i = X'_i$ and the result follows.
 
 (ii)--(iv) These equations can be verified  in a similar fashion.
\hfill $\Box $  \\ 

\section{The linear transformations $B$ and $B^*$}

\begin{definition}
\label{def:b}
\label{def:k}
\rm
Assume $V$ is a finite-dimensional
irreducible 
$U^{\geq 0}$-module.
Let the decompositions
$W_0, \ldots, W_d$ and
$W^*_0, \ldots, W^*_d$ be as in
Definition \ref{def:wdec}. 
\begin{enumerate}
\item 
Let $B:V\rightarrow V$ denote the 
linear transformation
such that for $0 \leq i \leq d$, $W_i$  
is an eigenspace of $B$ with eigenvalue 
$q^{2i-d}$.
\item
Let $B^*:V\rightarrow V$ denote the
linear transformation
such that for $ 0 \leq i \leq d$,  $W^*_i$ is  
 an eigenspace of $B^*$ with eigenvalue 
$q^{d-2i}$. 
\end{enumerate}
\end{definition}

Next we show that  $A,A^*,
B,B^*$ satisfy $q$-Weyl relations.   
\begin{lemma}
\label{thm:key}  
Assume $V$ is a finite-dimensional irreducible
$U^{\geq 0}$-module of type $\alpha$.
For
 $A,A^*$ from
Definition
\ref{def:aas}
and  $B,B^*$ from
Definition
\ref{def:b},
\begin{eqnarray}
\frac{qAB-q^{-1}BA}{q-q^{-1}}&=&\alpha I,
\label{eq:eq1}
\\
\frac{qBA^*-q^{-1}A^*B}{q-q^{-1}}&=&\alpha^{-1}I,
\label{eq:eq2}
\\
\frac{qA^*B^*-q^{-1}B^*A^*}{q-q^{-1}}&=&\alpha^{-1}I,
\label{eq:eq4}
\\
\frac{qB^*A-q^{-1}AB^*}{q-q^{-1}}&=&\alpha I.
\label{eq:eq3}
\end{eqnarray}
\end{lemma} 
\noindent {\it Proof:}
Let the decomposition $W_0, \ldots, W_d$ be as in
Lemma 
\ref{lem:wone}. 
To prove
(\ref{eq:eq1}),  
we show that 
$qAB-q^{-1}BA-\alpha (q-q^{-1})I$ 
is 0 on  $W_i$ for
$0 \leq i \leq d$. 
Since $B-q^{2i-d}I$ vanishes on
 $W_i$
by Definition
\ref{def:b}(i), 
  so does 
\begin{eqnarray}
\label{eq:p1}
(A-\alpha q^{d-2i-2}I)(B-q^{2i-d}I).
\end{eqnarray}
By
Lemma
\ref{lem:aw} (i), 
\begin{eqnarray}
\label{eq:p2}
(B-q^{2i+2-d}I)(A-\alpha q^{d-2i}I)
\end{eqnarray}
vanishes on $W_i$.
Subtracting $q^{-1}$times 
(\ref{eq:p2}) from
$q$ times (\ref{eq:p1}) 
we find that 
$qAB-q^{-1}BA-\alpha (q-q^{-1})I$ is 0 on
 $W_i$.   Relation 
 (\ref{eq:eq1}) follows.  
Lines 
(\ref{eq:eq2})--(\ref{eq:eq3})
can be proved in a similar manner.
\hfill $\Box $ \\

\section{The action of $B$ and $B^*$ on the decompositions}

\noindent In this section we describe how the maps
$B,B^*$ act on
our five decompositions. 

\begin{lemma}
\label{lem:baction}
Let $V$ denote a 
finite-dimensional irreducible $U^{\geq 0}$-module with 
weight space decomposition $U_0, \ldots, U_d$.
 Assume  that the decompositions  
 $V_0, \ldots, V_d$ and 
$V^*_0, \ldots, V^*_d$ are as 
in  Definition
\ref{def:vivis}. Let the decompositions 
$W_0, \ldots, W_d$
and $W^*_0, \ldots, W^*_d$  be as 
in Definition
\ref{def:wdec}.
Assume
$B,B^*$ are as in
Definition
\ref{def:b}.  Then the following hold
for $0 \leq i \leq d$. 
\begin{enumerate}
\item[{\rm (i)}] 
	$ (B-q^{d-2i}I)V_i\subseteq V_{i-1}$  and 
	$(B-q^{2i-d}I)V^*_i \subseteq V^*_{i-1}$.
\item[{\rm (ii)}] 
	$(B^*-q^{d-2i}I)V_i \subseteq V_{i+1}$ and $ (B^*-q^{2i-d}I)V^*_i\subseteq V^*_{i+1}$.
\item[{\rm (iii)}] 
        $(B-q^{2i-d}I)U_i\subseteq U_{i-1}$.
\item[{\rm (iv)}] 
	$(B^*-q^{d-2i}I)U_i \subseteq U_{i+1}$.
\item[{\rm (v)}] 
	$BW^*_i\subseteq W^*_{i-1}+W^*_i + W^*_{i+1}$.
\item[{\rm (vi)}] 
	$B^*W_i \subseteq W_{i-1}+W_i + W_{i+1}$.
\end{enumerate}
\end{lemma}

\noindent {\it Proof:} 
(i), (ii) Let $\alpha $ denote the type of $V$.
Recall that for  $0 \leq i \leq d$,  
$V_i$ is the eigenspace for $A$ with
eigenvalue $\alpha q^{2i-d}$.
To obtain the first half of (i) 
it is sufficient to show that 
\begin{eqnarray}
(A-\alpha q^{2i-2-d}I)(B-q^{d-2i}I)
\label{eq:want}
\end{eqnarray}
vanishes on $V_i$ for $0 \leq i \leq d$.
Since   $A-\alpha q^{2i-d}I$ vanishes on
$V_i$, so does 
\begin{eqnarray}
(B-q^{d-2i+2}I)(A-\alpha q^{2i-d}I).
\label{eq:ba}
\end{eqnarray}  Equation 
(\ref{eq:eq1}) implies  that 
\begin{eqnarray}
qAB-q^{-1}BA-\alpha (q-q^{-1})I
\label{eq:qab}
\end{eqnarray}
is 0 on $V_i$.
Adding
(\ref{eq:ba})
to 
$q$ times 
(\ref{eq:qab}),  
we find that
(\ref{eq:want}) vanishes on $V_i$.
Consequently,  $(B-q^{d-2i}I)V_i \subseteq V_{i-1}$ 
for  $0 \leq i \leq d$.

The second half of (i) and the relations in (ii) can be established similarly.  

(iii) We have
$$ \begin{array}{cccc} 
(B-q^{2i-d}I)U_i  &\subseteq &
(B-q^{2i-d}I)(U_0+\cdots +U_{i}) \hfill & \\
&=& (B-q^{2i-d}I)(V^*_0+\cdots + V^*_i)\hfill 
&\quad  (\mbox{by
Lemma \ref{lem:uvsv}(ii)})
\\
&\subseteq& V^*_0+\cdots +V^*_{i-1}\hfill &
\quad  (\mbox{by
Lemma \ref{lem:baction}(i)}
 )
\\
&=& U_0+\cdots +U_{i-1}\hfill&
\quad (\mbox{by
Lemma \ref{lem:uvsv}(ii)}
),
\end{array}$$
and also
$$ \begin{array}{cccc} 
(B-q^{2i-d}I)U_i  &\subseteq &
(B-q^{2i-d}I)(U_i+\cdots + U_{d})\hfill & 
\\
&=& (B-q^{2i-d}I)(V_i+\cdots +V_d)\hfill &
\quad \   (\mbox{by
Lemma \ref{lem:uvsv}(i)}
)
\\
&\subseteq& V_{i-1}+\cdots +V_d \hfill &
\quad \ (\mbox{by
Lemma \ref{lem:baction}(i)}
 )
\\
&=& U_{i-1}+\cdots +U_d\hfill &
\quad \ (\mbox{by
Lemma \ref{lem:uvsv}(i)}
).\end{array}$$
Combining these observations  
we find
$(B-q^{2i-d}I)U_i \subseteq U_{i-1}$.
 
(iv) To obtain this part, imitate the argument for (iii).   

(v) We have 
$$ \begin{array}{cccc} 
B W^*_i &\subseteq&
B (W^*_0+ \cdots + W^*_i) \hfill& 
\\
&=& B (V^*_{d-i}+ \cdots + V^*_d)\hfill&\hskip 1 truein  \quad (\mbox{by 
Lemma
\ref{lem:wsum}(iv)}
)
\\
&\subseteq& V^*_{d-i-1}+ \cdots + V^*_d \hfill&\hskip 1 truein \quad (\mbox{by 
Lemma \ref{lem:baction}(i)}
)
\\
&=& 
W^*_0+ \cdots + W^*_{i+1}\hfill&\hskip 1 truein \quad (\mbox{by 
Lemma
\ref{lem:wsum}(iv)}
), 
\end{array}$$
and also
$$ \begin{array}{cccc} 
B W^*_i &\subseteq&
B (W^*_i+ \cdots + W^*_d)\hfill&
\\
&=&B (V_i+ \cdots + V_d)\hfill&\hskip 1 truein
\quad (\mbox{by 
Lemma
\ref{lem:wsum}(iii)}
)
\\
&\subseteq& V_{i-1}+ \cdots + V_d\hfill&\hskip 1 truein\quad (\mbox{by 
Lemma \ref{lem:baction}(i)}
)
\\
&=& 
W^*_{i-1}+ \cdots + W^*_d\hfill&\hskip 1 truein\quad (\mbox{by 
Lemma
\ref{lem:wsum}(iii)}
).
\end{array}$$
Together these relations imply 
$B W^*_i \subseteq 
 W^*_{i-1} 
+
 W^*_i 
+
 W^*_{i+1}$.
  
 (vi) This argument is identical to (v).   \hfill $\Box $ \\

\section{Some relations involving $B,B^*, K^{\pm 1}$}

\noindent In this section we show
$B,B^*, K^{\pm 1}$ satisfy $q$-Weyl relations. 
\begin{lemma}
\label{thm:kkey}
Let $V$ be a finite-dimensional irreducible
$U^{\geq 0}$-module of type $\alpha $. 
Assume $B,B^*$ are as in
Definition
\ref{def:b}. Then both  
\begin{eqnarray}
\frac{qBK^{-1}-q^{-1}K^{-1}B}{q-q^{-1}}&=&\alpha^{-1} I,  
\label{eq:keq2}
\\
\frac{qB^*K-q^{-1}KB^*}{q-q^{-1}}&=&\alpha I.
\label{eq:keq3}
\end{eqnarray}
\end{lemma}
\noindent {\it Proof:}
Let $U_0, \ldots, U_d$ denote the weight space decomposition
for $V$.
Recall that for $0 \leq i \leq d$, $U_i$ is an eigenspace for
$K$ with eigenvalue $\alpha q^{2i-d}$.
To obtain (\ref{eq:keq2}),
we show
$qBK^{-1}-q^{-1}K^{-1}B-\alpha^{-1}(q-q^{-1})I$ vanishes on $U_i$ for
$0 \leq i \leq d$. 
Observe that $K^{-1}-\alpha^{-1}q^{d-2i}I$ vanishes on
 $U_i$
so
\begin{eqnarray}
\label{eq:kp1n}
(B-q^{2i-d-2}I)(K^{-1}-\alpha^{-1}q^{d-2i}I)
\end{eqnarray}
is 0 on $U_i$. From Lemma \ref{lem:baction} (iii) we see that 
$(B-q^{2i-d}I)U_i \subseteq U_{i-1}$. Therefore 
\begin{eqnarray}
\label{eq:kp2n}
(K^{-1}-\alpha^{-1}q^{d-2i+2}I)(B-q^{2i-d}I)
\end{eqnarray}
vanishes on $U_i$.
Subtracting $q^{-1}$times 
(\ref{eq:kp2n}) from
$q$ times (\ref{eq:kp1n}) 
we find
$qBK^{-1}-q^{-1}K^{-1}B-\alpha^{-1}(q-q^{-1})I$ vanishes on $U_i$.
Equation 
(\ref{eq:keq2}) follows, 
and relation 
(\ref{eq:keq3}) is proved similarly.
\hfill $\Box $ \\

\section{The $q$-Serre relations}

\noindent Next we 
show that the elements $B,B^*$ from Definition
\ref{def:b}
satisfy the $q$-Serre relations.

\begin{theorem}
\label{thm:qsab}
Let $V$ be a finite-dimensional irreducible $U^{\geq 0}$-module. 
Then the transformations  $B,B^*$  in Definition
\ref{def:b} satisfy the relations
\begin{eqnarray}
B^3B^*-\lbrack 3 \rbrack  B^2B^*B+\lbrack 3 \rbrack BB^*B^2-B^*B^3&=&0,
\label{eq:bserre}
\\
B^{*3}B-\lbrack 3 \rbrack B^{*2}BB^*+\lbrack 3 \rbrack B^*BB^{*2}-BB^{*3}&=&0.
\label{eq:bsserre}
\end{eqnarray}
\end{theorem}
\noindent {\it Proof:}
Let the decomposition $W_0, \ldots, W_d$
be as in 
Definition
\ref{def:wdec}.
Recall that for $0 \leq i \leq d$, 
$W_i$
 is an eigenspace for $B$
with eigenvalue $q^{2i-d}$.
In order to prove
\eqref{eq:bserre} we show that the
transformation  
$\Psi := 
B^3B^*-\lbrack 3\rbrack B^2B^*B+\lbrack 3\rbrack BB^*B^2-B^*B^3$ is
0 on  $W_i$  for $0 \leq i \leq d$.
Let $i$ be given and pick $v \in W_i$.
Observe that $B^*v \in 
W_{i-1}+
W_{i}+
W_{i+1}$ by
Lemma
\ref{lem:baction}(vi).
Next note that  
$(B-q^{2i-2-d}I)W_{i-1}=0$,
$(B-q^{2i-d}I)W_{i}=0$, and
$(B-q^{2i+2-d}I)W_{i+1}=0$.
By these comments
\begin{eqnarray*}
(B-q^{2i-2-d}I)
(B-q^{2i-d}I)
(B-q^{2i+2-d}I)
 B^*v=0.
\end{eqnarray*}
Therefore, 
\begin{eqnarray*}
\Psi v &=&
\bigl(B^3B^*-\lbrack 3\rbrack B^2B^*B+\lbrack 3\rbrack BB^*B^2-B^*B^3\bigr)v
\\
&=&
\bigl(B^3B^*-\lbrack 3\rbrack B^2B^*q^{2i-d}+
\lbrack 3\rbrack BB^*q^{4i-2d}-B^*q^{6i-3d}\bigr)v
\\
&=& 
\bigl(B-q^{2i-2-d}I\bigr)
\bigl(B-q^{2i-d}I\bigr)
\bigl(B-q^{2i+2-d}I\bigr)B^*v
\\
&=& 0.
\end{eqnarray*}
We have now shown $
\Psi W_i=0$ for $0 \leq i \leq d$.
Consequently
 $\Psi =0$ and 
(\ref{eq:bserre}) follows.
Line (\ref{eq:bsserre}) can be proved similarly.
\hfill $\Box $ \\

\section{
The proof of Theorem 
\ref{thm:nnmain}}

This section is devoted to a proof of 
Theorem
\ref{thm:nnmain}.

\begin{definition}
\label{def:rl}
Let $V$ be a finite-dimensional irreducible $U^{\geq 0}$-module
of type
$\alpha$.
Using the transformations $B,B^*$ in 
Definition
\ref{def:b}, 
 we introduce linear transformations 
$r:V\rightarrow V$ and $l:V\rightarrow V$ as follows:
\begin{eqnarray*}
r&=& \frac{\alpha I- KB^*}{q(q-q^{-1})^2},
\label{eq:rdef}
\\
l&=& \frac{\alpha^{-1}I- K^{-1}B}{q(q-q^{-1})^2}.
\label{eq:ldef}
\end{eqnarray*}
\end{definition}

\begin{lemma}
\label{lem:bKr} 
With reference to Definition \ref{def:rl}, 
we have
\begin{eqnarray*}
\label{eq:bsolve}
B&=&\alpha^{-1} K-q(q-q^{-1})^2Kl,
\\
\label{eq:bssolve}
B^*&=&\alpha K^{-1}-q(q-q^{-1})^2K^{-1}r.
\end{eqnarray*}
\end{lemma}  

\medskip 

\begin{theorem}
\label{lem:uqmoreaction}
Let $V$ denote a finite-dimensional irreducible
$U^{\geq 0}$-module of type $\alpha$.
Then the generators $R,L,K^{\pm 1}$ of 
$U^{\geq 0}$, together with
$r,l$ from
 Definition
\ref{def:rl},  satisfy the following relations on $V$:
\begin{eqnarray}
\label{eq:kkin2}
&& KK^{-1} = K^{-1}K = 1,\\
&&K R K^{-1} = q^2 R,
\hskip 1 truein 
K L K^{-1} = q^{-2} L,
\label{eq:app1}
\\
&&
K r K^{-1}= q^2 r , \hskip 1.1 truein  
 K l K^{-1}= q^{-2} l ,  
\label{eq:appr}
\\
&&
rL-Lr = {{\alpha^{-1}K-\alpha K^{-1}}\over {q-q^{-1}}},
\quad \quad 
lR-Rl = {{\alpha K^{-1}-\alpha^{-1}K}\over {q-q^{-1}}},
\qquad 
\label{eq:app3}
\\
&&
lL = Ll , \hskip 1.55 truein  rR=Rr,
\label{eq:app4}
\\
&&
0= R^3L - \lbrack 3 \rbrack R^2LR+ \lbrack 3 \rbrack RLR^2 - LR^3, 
\label{eq:app5}
\\
&&
0= L^3R - \lbrack 3 \rbrack L^2RL+ \lbrack 3 \rbrack LRL^2 - RL^3.
\label{eq:app6}
\\
&&
0= r^3l - \lbrack 3 \rbrack r^2lr+ \lbrack 3 \rbrack rlr^2 - lr^3, 
\label{eq:app7}
\\
&&
0= l^3r - \lbrack 3 \rbrack l^2rl+ \lbrack 3 \rbrack lrl^2 - rl^3.
\label{eq:app8}
\end{eqnarray}
\end{theorem}

\noindent {\it Proof:}
The relations in (\ref{eq:kkin2}), 
(\ref{eq:app1}) 
are defining relations 
(\ref{eq:nnbuq1})--(\ref{eq:nnbuq4}) of $U^{\geq 0}$.
To obtain
(\ref{eq:appr}), evaluate each of
(\ref{eq:keq2}), (\ref{eq:keq3})
using
Lemma \ref{lem:bKr}.
{F}or
(\ref{eq:app3}),
(\ref{eq:app4}), use Definition
\ref{def:aas},  Lemma \ref{lem:bKr},  and equations 
(\ref{eq:app1}) and 
(\ref{eq:appr}) to
evaluate 
(\ref{eq:eq1})--(\ref{eq:eq3}). 
Lines 
(\ref{eq:app5}),
(\ref{eq:app6}) are just defining relations  
(\ref{eq:nnbuq5}),  
(\ref{eq:nnbuq7}) respectively.
Finally, to demonstrate 
(\ref{eq:app7}) and 
(\ref{eq:app8}),
substitute the expressions in 
Lemma \ref{lem:bKr} into (\ref{eq:bserre}), (\ref{eq:bsserre}),  and 
apply relations 
(\ref{eq:app1}),
(\ref{eq:appr}).
\hfill $\Box $  
 
\begin{theorem}
\label{thm:modaction}
Let $V$ be a finite-dimensional irreducible 
$U^{\geq 0}$-module of type $\alpha $.
Assume the maps $r,l$ are as in Definition
\ref{def:rl}, and 
let $\varepsilon_0, \varepsilon_1$ denote scalars in $\lbrace 1,-1\rbrace $.
Then
$V$ supports an irreducible
$\uqhat$-module structure of  type $(\varepsilon_0, \varepsilon_1)$
for which the Chevalley generators act as follows:

\medskip
\centerline{
\begin{tabular}[t]{c|cccccccc}
        {\rm generator}  
       	& $e_0^+$  
         & $e_1^+$  
         & $e_0^-$  
         & $e_1^-$  
         & $K_0$  
         & $K_1$  
         & $K_0^{-1}$  
         & $K_1^{-1}$  
	\\
	\hline 
{\rm action on $V$} 
&  $R$ & $L$ & $\varepsilon_0 l$ & $\varepsilon_1 r$ &
$\varepsilon_0 \alpha^{-1}K$ & $\varepsilon_1 \alpha K^{-1}$ &
$\varepsilon_0\alpha K^{-1}$ 
& $\varepsilon_1\alpha^{-1}K$ 
\end{tabular}}
\end{theorem}
\noindent {\it Proof:}
To see that the above action 
on $V$ determines a
$\uqhat$-module,
compare equations 
(\ref{eq:kkin2})--(\ref{eq:app8}) 
with the defining relations for 
$\uqhat$ 
in
Definition \ref{def:uq}.
The
$\uqhat$-module
$V$ is irreducible, 
since $V$ is an irreducible $U^{\geq 0}$-module.
It is straightforward to check that $V$ has
 type
 $(\varepsilon_0, \varepsilon_1)$.
\hfill $\Box $  \\

\noindent {\bf Proof of Theorem \ref{thm:nnmain}}:
The ``existence'' part is immediate from
Theorem
\ref{thm:modaction}.
Concerning the ``uniqueness''  assertion,
we assume there exists
a $\uqhat$-module structure on
$V$ such that
the transformations
$e^+_0-R$, \ $e^+_1-L$, $K^{\pm 1}_0-
\varepsilon_0\alpha^{\mp 1}K^{\pm 1}$, and 
 $K^{\pm 1}_1
-\varepsilon_1\alpha^{\pm 1} K^{\mp 1}$ vanish on $V$.
Observe that 
$e^+_0$, $e^+_1$, $K^{\pm 1}_0$,
$ K^{\pm 1}_1$
 act on  $V$ according to the
table of
Theorem \ref{thm:modaction}. 
We show 
the remaining Chevalley generators 
$e^-_0, e^-_1$ must act on $V$
according to that table also.

First note that 
the given
$\uqhat$-module structure is irreducible,
since 
$R,L,K^{\pm 1}$ is a generating set for
$U^{\geq 0}$ and  
$V$ is irreducible as a $U^{\geq 0}$-module.
Next observe that  $K_0$ acts on $V$ as $\varepsilon_0 \alpha^{-1} K$.
Comparing 
(\ref{eq:nnkmove}) and 
(\ref{eq:s1dec}),  
we find that the weight space decomposition of the $U^{\geq 0}$-module
$V$ relative to $K$ coincides with the 
weight space decomposition of the
$\uqhat$-module $V$ relative to $K_0$, $K_1$.
Let $U_0, \ldots, U_d$ denote this (common) 
weight space decomposition.
Assume $l,r$ are as in Definition
\ref{def:rl}.
To show $e_0^-$ acts on $V$ as $\varepsilon_0 l$,
we set $W=\lbrace v \in V \mid (e_0^--\varepsilon_0 l)v=0 \rbrace $, and 
argue that  $W=V$.
For this,  it suffices to prove that $W\not=0$, and $W$
is invariant under the operators 
$R,L, K^{\pm 1}$.
Using the right-hand equation in
(\ref{eq:appr}),  we find that 
$l U_i \subseteq U_{i-1}$ for $0 \leq i\leq d$.
In particular,  
$lU_0=0$.
Also
$e^-_0U_0=0$ by
(\ref{eq:altemove2})
so
$e_0^- -\varepsilon_0 l$ vanishes on $U_0$. Therefore 
$U_0 \subseteq W$,  so
$W\not=0$.
By (\ref{eq:buq5}) (with $i=0$),
the second relation  in (\ref{eq:app3}),
and the fact that both $e_0^+-R$ and  $K_0-\varepsilon_0 \alpha^{-1} K$ vanish on $V$, 
we deduce that the commutator
$\lbrack e_0^- -\varepsilon_0 l,R\rbrack$ vanishes on $V$.
This implies 
$RW\subseteq W$. 
By 
(\ref{eq:buq6}) we have
$\lbrack e_0^-, e^+_1 \rbrack =0$.
Therefore, 
the first equation in (\ref{eq:app4}),
combined with the fact that $e_1^+-L$ vanishes on $V$
implies that
$\lbrack e_0^- -\varepsilon_0 l,L\rbrack$ vanishes on $V$.
Using this,  we find that 
$LW\subseteq W$.
Observe
$K_0 e^-_0=q^{-2}e^-_0K_0$ by
(\ref{eq:buq3})
and  $K_0- \varepsilon_0 \alpha^{-1} K$ vanishes on $V$, 
so 
$K e^-_0-q^{-2}e^-_0K$ is 0 on $V$.
Combining this with the second equation in
(\ref{eq:appr})
we determine that 
$K(e^-_0 -\varepsilon_0 l)$ and
$q^{-2}(e^-_0 -\varepsilon_0 l)K$
agree on $V$.
Thus, 
$KW\subseteq W$ and then
$K^{-1}W\subseteq W$.
We have now shown that
$W\not=0$ and $W$ is invariant under each of
$R, L,
K^{\pm 1}$. Therefore 
$W=V$ since
$V$ is irreducible as a $U^{\geq 0}$-module.  
We conclude $(e_0^- -\varepsilon_0 l)V=0$ so $e_0^-$ acts on $V$ 
as $\varepsilon_0 l$.
By a similar argument we find  
$e_1^-$ acts on $V$ as $\varepsilon_1 r$.
 Consequently,  $e^-_0, e^-_1$ must act on
$V$ 
according to the table of
Theorem \ref{thm:modaction}. 
Hence,  the given
$\uqhat$-module structure 
is unique. We already showed
that this
$\uqhat$-module structure is irreducible and it
clearly has type
$(\varepsilon_0, \varepsilon_1)$.

\hfill $\Box $  \\ 

\section{The proof of Theorem
\ref{thm:nnmain2}}

This section is devoted to a proof of
 Theorem 
\ref{thm:nnmain2}.
We begin with a few comments about
the quantum algebra $U_q(\mathfrak{sl}_2)$ and its modules.

\begin{definition}
\rm \cite[p.~122]{Kassel}
\label{def:udef}
The quantum algebra $U_q(\mathfrak{sl}_2)$ is the unital associative
$\K$-algebra with generators $e^{\pm }$, $k^{\pm 1}$
which satisfy the following relations:
\begin{eqnarray*}
kk^{-1}&=&k^{-1}k = 1,
\\
ke^{\pm }k^{-1} &=& q^{\pm 2} e^{\pm },
\\
\lbrack e^+, e^-\rbrack &=& \frac{k-k^{-1}}{q-q^{-1}}.
\end{eqnarray*}
\end{definition}

\begin{lemma}
\rm \cite[p.~128]{Kassel}
\label{lem:fdsl} 
If  $V$ is a finite-dimensional irreducible 
$U_q(\mathfrak{sl}_2)$-module, then  there exist $\varepsilon \in \lbrace 1,-1\rbrace$
and a basis 
$v_0, v_1, \ldots, v_d$ for $V$ 
such that $kv_i = \varepsilon q^{2i-d}v_i$
for $0 \leq i \leq d$, \
$e^+v_i = \lbrack i+1\rbrack v_{i+1}$
for $0 \leq i \leq d-1$, 
$e^+v_d =0$,\
$e^-v_i = \varepsilon\lbrack d-i+1\rbrack v_{i-1}$
for $1 \leq i \leq d$, and 
$e^-v_0 =0$.
\end{lemma}

The proof of the next lemma is straightforward. 

\begin{lemma}
\label{lem:uqgives}
Let $V$ denote a finite-dimensional irreducible
$\uqhat$-module.
Then for $i\in \lbrace 0,1\rbrace $,
there exists a unique
$U_q(\mathfrak{sl}_2)$-module structure on $V$ such that
$e^{\pm } -e^{\pm}_i$ and 
$k^{\pm 1}-K^{\pm 1}_i$
vanish on $V$.
\end{lemma}

In our proof of Theorem
\ref{thm:nnmain2}, we will use the following facts 
concerning
$\uqhat $-modules.

\begin{lemma} 
\label{lem:tp}
Let $V$ be a finite-dimensional irreducible
$\uqhat$-module,  and let $U_0, \ldots, U_d$
denote the corresponding weight space decomposition relative
to $K_0$, $K_1$.
Then the following hold  for $0 \leq j \leq d/2$.
\begin{enumerate}
\item[{\rm (i)}]  The restriction
of $(e^+_0)^{d-2j}$ to $U_j$ is an isomorphism of
vector spaces from $U_j $ to $U_{d-j}$.
\item[{\rm (ii)}] 
For all $v \in U_j$,  \ $(e^+_0)^{d-2j+1}v=0$ if
and only if $e^-_0v=0$.
\item[{\rm (iii)}] 
The restriction
of $(e^+_1)^{d-2j}$ to $U_{d-j}$ is an isomorphism of
vector spaces from $U_{d-j} $ to $U_j$.
\item[{\rm (iv)}] 
For all $v \in U_{d-j}$, \ $(e^+_1)^{d-2j+1}v=0$ if
and only if $e^-_1v=0$.
\end{enumerate}
\end{lemma}
\noindent {\it Proof:}
For (i) and (ii),  we view $V$ as 
a 
$U_q(\mathfrak{sl}_2)$-module via Lemma
\ref{lem:uqgives} (with $i=0$).  
As a $U_q(\mathfrak{sl}_2)$-module, 
$V$ is  a direct sum of irreducible $U_q(\mathfrak{sl}_2)$-modules (see
for example, \cite[p. 144]{Kassel}). 
Let $W$ denote one of the irreducible 
$U_q(\mathfrak{sl}_2)$-module summands.
Applying
Lemma
\ref{lem:fdsl}
to $W$, we find there exists an integer $r$ $(0 \leq r \leq d/2)$
such that for $0 \leq i \leq d$,
$W\cap U_i$ is one-dimensional if  
$r \leq i \leq d-r$ and is zero otherwise.
Moreover
$e^+_0 (W\cap U_i) = W\cap U_{i+1}$ for $r \leq i \leq d-r-1$, \
$e^+_0 (W\cap U_{d-r})=0$,  \ 
$e^-_0 (W\cap U_i) = W\cap U_{i-1}$ for $r+1 \leq i \leq d-r$, \ and 
$e^-_0 (W\cap U_r)=0$.
Results (i), (ii)  follow.   The statements in (iii) and (iv) can be shown 
in exactly the same way. 
\hfill $\Box $ \\
  
\noindent {\bf  Proof of Theorem 
\ref{thm:nnmain2}:} Let $V$ be a finite-dimensional irreducible $\uqhat$-module
of type $(\varepsilon_0, \varepsilon_1)$.  
We first prove that the desired $U^{\geq 0}$-module structure on $V$ 
exists. 
Let $R,L,K^{\pm 1}$ act on $V$ as
$e^+_0$, $e^+_1$, $\varepsilon_0 \alpha^{\pm 1}K^{\pm 1}_0$
respectively.  Then 
using the defining relations for
$\uqhat$ in 
Definition
\ref{def:uq},
it is easy to see that 
$R,L,K^{\pm 1}$ satisfy
(\ref{eq:nnbuq1})--(\ref{eq:nnbuq7}),  and therefore
induce a $U^{\geq 0}$-module structure on $V$.
{F}rom the construction,  the transformations 
$e^+_0-R, e^+_1-L, K^{\pm 1}_0-
\varepsilon_0
\alpha^{\mp 1}K^{\pm 1}$ vanish on $V$.
Since  $K^{\pm 1}_0-
\varepsilon_0
\alpha^{\mp 1}K^{\pm 1}$ vanish on $V$,
and since
$K_0K_1-\varepsilon_0 \varepsilon_1I$ also vanishes
on $V$  by
Lemma \ref{lem:nnweight},
we find that
 $K^{\pm 1}_1
 -\varepsilon_1
 \alpha^{\pm 1} K^{\mp 1}$ vanish on $V$. 
We have now shown the desired
 $U^{\geq 0}$-module structure exists,  and
 it is clear this
 $U^{\geq 0}$-module structure is unique.
 Next we show the
 $U^{\geq 0}$-module structure
  is irreducible.
 Let $W$ denote an irreducible 
$U^{\geq 0}$-submodule of  
$V$.  Then 
$W$ is invariant under  the actions of 
$e^+_0$, $e^+_1$, $K^{\pm 1}_0$,  and  
 $K^{\pm 1}_1$.
We argue that $e^-_0W \subseteq W$
and 
$e^-_1W \subseteq W$.
To demonstrate the first of these assertions, it will be enough to show that 
${\widetilde W} = \lbrace w \in W \mid e^-_0w \in W\rbrace$ is nonzero
and
 invariant
under each of  $R,L, K^{\pm 1}$.  
It follows {f}rom
(\ref{eq:buq5}) (with $i=0$) 
and the fact that  $e^+_0-R$ and
$\varepsilon_0K^{\pm 1}_0-
\alpha^{\mp 1}K^{\pm 1}$ \
are 0 on $V$ that 
$R{\widetilde W}\subseteq {\widetilde W}$.
Also $\lbrack e^-_0, e^+_1\rbrack =0$ by (\ref{eq:buq6})   
and  $e^+_1-L$ vanishes on $V$, so that 
$L{\widetilde W}\subseteq {\widetilde W}$.
By 
(\ref{eq:buq3})  we have
$K_0e^{-}_0K^{-1}_0 = q^{-2}e^{-}_0$. 
By this and since
$K_0-
\varepsilon_0\alpha^{-1}K$ is 0 on $V$,
we find that
$Ke^{-}_0 - q^{-2}e^{-}_0K$ vanishes on $V$. 
Consequently, 
 $K{\widetilde W}\subseteq {\widetilde W}$, 
 and then
 $K^{-1}{\widetilde W}\subseteq {\widetilde W}$ holds as well.
 
To verify that  ${\widetilde W}\not=0$,
let $U_0, \ldots, U_d$ denote the
weight space decomposition for the 
$\uqhat$-module $V$ relative to $K_0$, $K_1$.  
As $W$ is invariant under $K_0$ and $K_0K_1$ acts as $\varepsilon_0\varepsilon_1 I$ on $V$, 
 it must be that  $W=\sum_{i=0}^d (W\cap U_i)$.
Since $W\not=0$, there exists an integer
$i$ $(0 \leq i \leq d)$ such that
$W\cap U_i \not=0 $.
Define $r=\mbox{min}\lbrace i\mid 0 \leq i \leq d, W\cap U_i\not=0\rbrace$
and 
$s=\mbox{max}\lbrace i\mid 0 \leq i \leq d, W\cap U_i\not=0\rbrace$.
We  prove that 
 $r+s=d$.
Suppose for the moment that $r+s<d$.  As 
$r\leq s$,  we must have $r<d/2$.
Then for any nonzero $v \in W \cap U_r$,
$(e^+_0)^{d-2r}v$ is contained in
$W\cap U_{d-r}$ which is 0 because $d-r > s$,  so  
$(e^+_0)^{d-2r}v=0$, contradicting
Lemma
\ref{lem:tp}(i).
Next assume $r+s>d$.    Then since $r\leq s$,  we have  $s>d/2$.
For any nonzero $v \in W \cap U_s$,
$(e^+_1)^{2s-d}v$ is contained in
$W\cap U_{d-s} = 0$,   so  
$(e^+_1)^{2s-d}v=0$,  which contradicts 
Lemma
\ref{lem:tp}(iii).
Thus, $r+s=d$ must hold, and {f}rom       
$r\leq s$,  we deduce that  $r\leq d/2$.
Let  $v$ denote a nonzero vector in $W\cap U_r$.
Observe $(e^+_0)^{d-2r+1}v$ is contained in $W\cap U_{d-r+1}$, 
and 
$W\cap U_{d-r+1}=0$, 
so 
 $(e^+_0)^{d-2r+1}v=0$.
Applying Lemma \ref{lem:tp} (ii) to
the $\uqhat$-module $V$, 
we obtain that 
$e^-_0v=0$.  Therefore, 
$v$ is a nonzero element of ${\widetilde W}$.
We have now shown $\widetilde W$ is nonzero and invariant under
each of the operators $R,L,K^{\pm 1}$.
Consequently,  by the irreducibility of $W$ as a $U^{\geq 0}$-module,
 ${\widetilde W}=W$ must hold.  
Therefore $e^-_0W \subseteq W$.
In just the same fashion, 
 $e^-_1W \subseteq W$, so that $W$ is a $\uqhat$-submodule of  $V$.
By  construction,  
$W\not=0$ 
so $W=V$.
We conclude that the 
$U^{\geq 0}$-module structure on $V$ is irreducible.
It is routine to show the 
$U^{\geq 0}$-module structure on $V$ has type
$\alpha$.
\hfill $\Box $ \\

\section{Irreducible $\uqhat^{\geq 0}$-modules} 

In this section we compare the finite-dimensional irreducible
$\uqhat$-modules with the  
finite-dimensional irreducible
$\uqhat^{\geq 0}$-modules.

Let  $V$ be a finite-dimensional irreducible module for $\uqhat^{\geq 0}$.
The central element  $K_0K_1$ must act as some scalar $\gamma$ times the identity
map on $V$.   Arguing as in Lemma \ref{lem:s1}, we see that there exists a nonzero
scalar $\alpha \in \K$ and a decomposition $U_0,\dots, U_d$ of $V$ such that 
$(K_0-\alpha q^{2i-d})U_i = 0$ for $0 \leq i \leq d$.  
Then $(K_1- \gamma \alpha^{-1}q^{d-2i})U_i = 0$ for $0 \leq i \leq d$,
so each of  the spaces
$U_0, \ldots, U_d$
is a common eigenspace for
$K_0, K_1$.     Setting $\beta = \gamma \alpha^{-1}$, we say 
$V$ has type $(\alpha, \beta)$ and diameter $d$.       

Assume $V$ is a finite-dimensional irreducible 
$\uqhat^{\geq 0}$-module  
 of type $(\alpha,\beta)$.
   Then $V$ remains irreducible when regarded as a module
for the subalgebra of $\uqhat^{\geq 0}$ generated by $K_0^{\pm 1}$, $e_0^+$ and $e_1^+$.
Thus, $V$ admits the structure of an irreducible $U^{\geq 0}$-module of type $\alpha$.   By Theorem
\ref{thm:nnmain},
for any choice of scalars $\varepsilon_0, \varepsilon_1$ {f}rom  $\{1,-1\}$,  
there exists a unique   
$\uqhat$-module structure on $V$ 
such that  the operators
$e^+_0-R, \ e^+_1-L, \ K^{\pm 1}_0-
\varepsilon_0\alpha^{\mp 1}K^{\pm 1}$, and
 $K^{\pm 1}_1
-\varepsilon_1\alpha^{\pm 1} K^{\mp 1}$ vanish  on $V$.
This $\uqhat$-module structure is irreducible and of
type
$(\varepsilon_0, \varepsilon_1)$.   When this $\uqhat$-module structure
on $V$ is then restricted to $\uqhat^{\geq 0}$, we will recover the original $\uqhat^{\geq 0}$-structure on $V$, provided
$\alpha = \varepsilon_0$ and $\beta = \varepsilon_1$.     

Next suppose that $V$ is a finite-dimensional
irreducible module for $\uqhat$.
We claim that $V$ remains irreducible as a module for
the subalgebra
$\uqhat^{\geq 0}$.
To see this,
let $W$ denote a nonzero
$\uqhat^{\geq 0}$-submodule of $V$.  
We show $W=V$.  By its definition, 
$W$ is invariant under each of 
$e_0^+, e_1^+, 
K_0^{\pm 1}, K_1^{\pm 1}$.
By Theorem \ref{thm:nnmain2}, for any nonzero  scalar $\alpha$ in $\K$, 
there is the structure of a  $U^{\geq 0}$-module
on $V$  such that the operators $e^+_0-R, \ e^+_1-L, \ K^{\pm 1}_0-
\varepsilon_0\alpha^{\mp 1}K^{\pm 1}$, and
 $K^{\pm 1}_1
-\varepsilon_1\alpha^{\pm 1} K^{\mp 1}$ vanish  on $V$.  
{F}rom this we see that 
$W$ is invariant under each of $R,L,K^{\pm 1}$,
and is therefore a 
 $U^{\geq 0}$-submodule of $V$.
But the $U^{\geq 0}$-module
structure on $V$ is irreducible by
Theorem
\ref{thm:nnmain2},
so $W=V$ and our claim is proved.
Note that if $V$ has type $(\varepsilon_0, \varepsilon_1)$ as a module for $\uqhat$,
then  $V$ has type $(\varepsilon_0, \varepsilon_1)$ as a module for $\uqhat^{\geq 0}$.  

Let us summarize these findings in our final result.

\begin{theorem}\label{thm:final} 
For any scalars $\varepsilon_0, \varepsilon_1$ taken from
the set $\lbrace 1,-1\rbrace$ the following hold.
\begin{enumerate}
\item[{\rm (i)}]   Let $V$ be a 
finite-dimensional irreducible
$\uqhat^{\geq 0}$-module
of type
$(\varepsilon_0,\varepsilon_1)$.
Then the action of 
 $\uqhat^{\geq 0}$  on $V$ extends uniquely
to an action of 
 $\uqhat$ on $V$.
The resulting 
$\uqhat$-module structure on $V$ is irreducible and
of type  $(\varepsilon_0,\varepsilon_1)$.
\item[{\rm (ii)}]   Let $V$ be a finite-dimensional irreducible $\uqhat$-module of type
$(\varepsilon_0,\varepsilon_1)$.
When the 
$\uqhat$-action is restricted to
 $\uqhat^{\geq 0}$,
the resulting $\uqhat^{\geq 0}$-module structure on $V$
is irreducible and of type
$(\varepsilon_0,\varepsilon_1)$.
\end{enumerate} \end{theorem}
     
\section{Acknowledgments}
The second author
would like to thank
Kenichiro Tanabe and Tatsuro Ito
for many conversations on
the subject of
$\uqhat$ and its modules.

\bigskip

\noindent Georgia Benkart  (email: {\tt benkart@math.wisc.edu})  \hfil\break
Paul Terwilliger (email:   {\tt terwilli@math.wisc.edu}) 
\medskip

\noindent Department of Mathematics \hfil\break
\noindent University of Wisconsin \hfil\break 
\noindent 480 Lincoln Drive \hfil\break
\noindent Madison, WI 53706-1388 USA 

\end{document}